\documentclass[a4paper,11 pt]{amsart}
\usepackage{srollens-en}[2006/10/02]
\usepackage{rotating, tabularx}

\title{The Kuranishi-space of complex parallelisable nilmanifolds}
\author{S\"onke Rollenske}
\address{S\"onke Rollenske\\Department of Mathematics\\
 Imperial College London \\
 SW7 2AZ London\\
 United Kingdom}
\email{s.rollenske@imperial.ac.uk}

\renewcommand{\lg}{\ensuremath{\gothg}}
\newcommand{\lh}{\gothh}
\newcommand{\einsnull}[1]{{{#1}^{1,0}}}
\newcommand{\nulleins}[1]{{{#1}^{0,1}}}

\newcommand{\Kur}{\mathrm{Kur}}
\newcommand{\obs}{\mathrm{obs}}
\newcommand{\refb}[1]{{\upshape (\ref{#1})}}

\newcommand{\margincom}[1]{\marginpar{\quad}}
\begin{document}

\begin{abstract}
We show that the deformation space of complex parallelisable nilmanifolds can be described by polynomial equations but is almost never smooth. This is remarkable since these manifolds have trivial canonical bundle and are holomorphic symplectic in even dimension. We describe the Kuranishi space in detail in several examples and also analyse when small deformations remain complex parallelisable
\end{abstract}

\maketitle
\tableofcontents

%
%
\section{Introduction}

Left-invariant geometric structures  on nilmanifolds, i.e., compact quotients of (real) nilpotent Lie groups, have proved to be both very rich and accessible for an in depth study. Thus many examples and counter-examples in (complex) differential geometry are of this type. 

In this paper we are concerned with deformations of complex structures for complex parallelisable nilmanifolds, which are the compact quotient of complex nilpotent Lie groups. \margincom{(explain more here or below??)}

The study of deformations of complex structures on compact complex manifolds has been an important topic since it was first developed by Kodaira and Spencer in \cite{Kod-sp58}. A deformation of a given compact complex manifold $X$ is a flat proper map $\pi:\kx\to \kb$ of (connected) complex spaces such that all the fibres are smooth manifolds together with an isomorphism with $X\isom \ky_0=\inverse\pi(0)$ for a point $0\in \kb$. If $\kb$ is a smooth $\pi$ is just a holomorphic submersion. Kodaira and Spencer showed that first order deformations correspond to elements in $H^1(X, \Theta_X)$ where $\Theta_X$ is the sheaf of holomorphic tangent vectors.

A key result is now the theorem of Kuranishi which, for a given compact complex manifold $X$, guarantees the existence of a locally complete space of deformations  $\kx\to\Kur(X)$ which is versal at the point corresponding to $X$. In other word, for every deformation $\ky\to \kb$ of $X$ there is a small neighbourhood $\ku$ of $0$ in $\kb$ yielding a diagram
\[\xymatrix{ \ky\restr{\ku}\isom f^*\kx \ar[d]\ar[r] & \kx \ar[d]\\
\ku\ar[r]^f&\Kur(X),}\]
and in addition the differential of $f$ at $0$ is unique.

The Kuranishi family $\Kur(X)$ hence parametrises all sufficiently small deformations of $X$. In general the map $f$ will not be unique which is roughly due to the existence of automorphisms.

Another point of view, which we will mainly adopt in this paper, is the following: consider $X$ as a differentiable manifold together with an integrable almost complex structure $(M,J)$, i.e., $J:TM\to TM$, $J^2=-\id_{TM}$ and the Nijenhuis integrability condition holds (see \refb{nijenhuis} below).

A deformation of $X$ can be viewed as a family of such complex structures $J_t$ depending on some parameter $t\in \kb$ with $J=J_0$. The construction of the Kuranishi space can then be made explicit after the choice of a hermitian metric on $M$. We will go through this construction for our special case of complex parallelisable nilmanifolds in Section \ref{kurani}.

In general the Kuranishi space can be arbitrarily bad but we can hope for better control over the deformations if we restrict our class of manifolds. If, for example,  $X$ is K\"ahler and has trivial canonical bundle, i.e., $X$ is a Calabi-Yau manifold, then the Tian-Todorov Lemma implies that the Kuranishi space is indeed smooth; we say that $X$ has unobstructed deformations. 
These manifolds are very important both in physics and in mathematics for example in the context of mirror symmetry.
This result fails if we drop the K\"ahler condition \cite{ghys95}.

The only nilmanifolds which can carry a K\"ahler structure are tori but it was proved by  Cavalcanti and Gualtieri \cite{caval-gual04} and, independently, by Babaris, Dotti and Verbitsky \cite{0712.3863v3} that nilmanifolds with left-invariant complex structure always have trivial canonical bundle. In addition, all known examples in the context of left-invariant complex structures on  nilmanifolds, e.g., complex tori, the Iwasawa manifold \cite{nakamura75}, Kodaira surfaces \cite{borcea84},  abelian complex structures  \cite{con-fin-poon06, mpps06} (see Section \ref{definitions} for a definition), had unobstructed deformations. Therefore it was speculated if this holds for all  left-invariant complex structures.
This was supported by results on weak homological mirror symmetry for nilmanifolds \cite{poon06,0708.3442v2}.

On the other hand Catanese and Frediani observed in their study of deformations of principal holomorphic torus bundles, which are in particular nilmanifolds with left-invariant complex structure, that the Kuranishi space can be singular \cite{cat-fred06}.

In this article we want to study the Kuranishi space of complex parallelisable nilmanifolds.
These were very intensively studied by Winkelmann \cite{winkelmann98} and they enjoy many interesting properties. We will only be concerned with their deformations.

We will show in particular, that the Kuranishi space of a complex parallelisable nilmanifold is  \emph{almost always} singular thus showing that no analog of the Tian-Todorv theorem can exist for nilmanifolds.

Nevertheless, the Kuranishi space can not become too ugly:
\begin{custom}[Theorem \ref{polynom}]
If $X=\Gamma\backslash G$ is a complex parallelisable nilmanifold and $G$ is $\nu$-step nilpotent, then $\Kur(X)$ is cut out by polynomial equations of degree at most $\nu$.
\end{custom}
In Section \ref{nonexam} we will give an example that the bound on the degree does not remain valid for general nilmanifolds but, as far as we know, there could be a larger bound depending on the step-length and the dimension only.

We believe that the Lie-algebra $\lg$ of $G$ cannot be too far from being free if the Kuranishi space is smooth and all examples that we found were actually free. Unfortunately, the analysis of the obstructions of higher order becomes very complicated but we can at least prove the following:
\begin{custom}[Theorem \ref{freecond}\slash Corollary \ref{b_m}]
Let $X=\Gamma\backslash G$ be a complex parallelisable nilmanifold and let $\lg$ be the Lie-algebra of $G$.
If  $\lg\slash [\lg,[\lg,\lg]]$ is not isomorphic to a free 2-step nilpotent Lie-algebra then there is a non-vanishing obstruction in degree 2 and the Kuranishi space is singular.

In particular, if $\lg$ is 2-step nilpotent then  $\Kur(X)$ is smooth if and only if  $\lg$ is a free 2-step nilpotent Lie-algebra.
\end{custom}

It is a natural question which infinitesimal deformations in $H^1(X,\Theta_X)$ integrate to a 1-parameter family of complex parallelisable complex structure and we show in Section \ref{remainparall} that this is the case if and only if they are infinitesimally complex parallelisable. The same results holds for abelian complex structures \cite{con-fin-poon06}.

From this we can also deduce that every complex parallelisable nilmanifold which is not a torus has small deformations which are no longer complex parallelisable (Corollary \ref{nondef}).
On the other hand it is known  that small deformations at least remain in the category of nilmanifolds with left-invariant complex structure (see Section \ref{kurani} or \cite{rollenske07b}).

In Section \ref{examples} we will give several explicit examples, mostly in small dimension. As far as we know, these are the first examples of compact complex manifolds with trivial canonical bundle (or even holomorphic symplectic structure) which have non-reduced Kuranishi-space.

\subsubsection*{Acknowledgements}
This research was carried out at Imperial College London supported by a DFG Forschungsstipendium\margincom{get precise formulation} and I would like to thank the Geometry group there for their hospitality. Fabrizio Catanese, Fritz Grunewald, Andrey Todorv and J\"org Winkelmann made several useful comments during a talk at the University of Bayreuth.

\margincom{Comments: Defined over $\IQ$. Non-degeneracy is the only condition.}

\section{Complex parallelisable nilmanifolds and nilmanifolds with left-invariant complex structure}\label{definitions}

Let $G$ be a  simply connected, complex, nilpotent Lie-group with Lie-algebra $\lg$ and $\Gamma\subset G$ a lattice, i.e., a discrete cocompact subgroup. By a theorem of Mal'cev \cite{malcev51} such a lattice exists if and only if the real Lie-algebra underlying $\lg$ can be defined over $\IQ$.

The most important invariant attached to a nilpotent Lie-algebra (or Lie-group) is its nilpotency index, also called step length. It is defined as follows: consider the descending central series, inductively defined by 
\[ \kc_0\lg:=\lg, \qquad \kc_{k+1}\lg=[\kc_k\lg,\lg].\]
Then $\lg$ is nilpotent if and only if there exists a $\nu$ such that $\kc^\nu\lg=0$. The smallest such $\nu$ is called the nilpotency index.

Since the multiplication in $G$ is holomorphic we can act with elements of $\Gamma$ on the left; the quotient $X:=\Gamma\backslash G$ is a complex parallelisable compact nilmanifold.

The nilpotent complex Lie-group $G$ acts transitively on $X$ by multiplication on the right and this is in fact an equivalent characterisation of $\IC$-parallelisable nilmanifolds \cite{wang54}.

As already remarked by Nakamura \cite{nakamura75} not all deformations of $\IC$-parallelisable nilmanifolds are again $\IC$-parallelisable but, as we will discuss in section \ref{kurani}, we can describe all deformations in the slightly more general framework of nilmanifolds with left-invariant complex structures which we will now explain.

Let $H$ be a simply connected, real, nilpotent Lie-group with Lie-algebra $\lh$  and containing a lattice $\Gamma$. Taking the quotient yields a real nilmanifold $M:=\Gamma\backslash H$.

An almost complex structure $J: \lh\to \lh $ defines an almost complex structure on $H$ by left-translation and this almost complex structure is integrable if and only if the Nijenhuis condition
\begin{equation}\label{nijenhuis}
 [x,y]-[Jx,Jy]+J[Jx,y]+J[x,Jy]=0
\end{equation}
holds for all $x,y\in \lh$. In this case we call the pair $(\lh, J)$ a Lie-algebra with complex structure.

The action of $\Gamma$ on the left is then holomorphic and we get an induced complex structure on $M$. We call  $(M,J)$  a nilmanifold with left-invariant complex structure.

Note that the multiplication in $H$ induces an action on the left on $M$ if and only if $\Gamma$ is normal if and only if $H=\IR^n$ is abelian; there is always an action on the right which is holomorphic if and only if $(H,J)$ is a complex Lie-group.

By abuse of notation we will call a tensor, e.g., a vector field, differential form or metric,  on $M$ left-invariant if its pullback to the universal cover $H$ is left-invariant.

The complexified Lie-algebra $\lh_\IC=\lh\tensor_\IR\IC$ decomposes as
\[ \lh_\IC=\einsnull\lh\oplus\nulleins\lh\]
where $\einsnull\lh$ is the $i$-eigenspace of $J$ and $\nulleins\lh=\overline{\einsnull\lh}$ is the $(-i)$-eigenspace.

It is not hard to see that the complex structure is integrable if and only if $\einsnull \lh$ is a (complex) Lie-subalgebra of $\lh_\IC$. 

The complex structure $J$ makes $(\lh,J)$ into a complex Lie-algebra if and only if the bracket is $J$-linear, i.e., for all $x,y\in \lh$ we have \begin{equation}\label{Clie}
[Jx,y]=J[x,y].
\end{equation}
In this case $H$ is a complex Lie-group and $(M,J)$ is $\IC$-parallelisable as above.
The following equivalent characterisation is also well known.

\begin{lem}\label{parallchar}
A Lie-algebra with  complex structure $(\lh,J)$ is a complex Lie-algebra if and only if $[\einsnull\lh, \nulleins\lh]=0$. 
In this case the canonical projection
\[\pi: (\lh,J)\to \einsnull\lh, \qquad z\mapsto \frac{1}{2}(z-iJz)\]
is an isomorphism of complex Lie algebras.
\end{lem}
\pf
Let $x,y\in \lh$ and consider $X:=\frac{1}{2}(x-iJx)\in \einsnull\lh$ and $\bar Y:=\frac{1}{2}(y+iJy)\in \nulleins\lh$. Then
\begin{align*}
  [X,\bar Y]&= \frac{1}{4}[x-iJx, y+iJy]\\&= \frac{1}{4}([x,y]-i^2[Jx,Jy]-i([Jx,y]-[x,Jy])\\&=\frac{1}{4}([x,y]+[Jx,Jy])-i([Jx,y]-[x,Jy])
\end{align*}
and we see that this vanishes if and only if
\[ [x,y]=-[Jx,Jy] \text{ and } [Jx,y]=[x,Jy].\]
If we combine these two equations with the Nijenhuis tensor \refb{nijenhuis} then we get the identity $-2[x,y]=2J[Jx,y]$ which becomes \refb{Clie} after applying $J$ to it and dividing by $-2$. On the other hand the equations are certainly fulfilled if \refb{Clie} holds and we have shown the claimed equivalence.

The second claim is proved by a similar computation: since $\pi$ is an isomorphism of complex vector spaces it remains to show that $\pi$ is a homomorphism of Lie-algebras. Indeed for $x,y\in\lh$ we have using \refb{Clie}
\[[\pi(x), \pi(y)]= \frac{1}{4}[x-iJx, y-iJy]= \frac{1}{4}([x,y]+i^2[Jx,Jy]-2iJ[x,y])=\pi([x,y]).\]
\qed

\begin{rem}[Notation]\label{notation}
In order to make our notation more transparent $\lh$, $H$ and $M$  will always denote a real Lie-algebra, Lie-group or nilmanifold, often equipped with a (left-invariant) complex structure $J$. We will only consider integrable complex structures.

The notations $\lg$, $G$ and $X$ will be reserved for their complex parallelisable counterparts. If we need to access the underlying real object with left-invariant complex structure we will write for example $\lg=(\lh,J)$. By the above Lemma we can then identify
\[ \lg_\IC=\lh_\IC=\lg\oplus \bar \lg\]
where the bracket on $\bar \lg$ is given by $[\bar x, \bar y]=\overline{[x,y]}$ and $[\lg, \bar\lg]=[\einsnull\lh, \nulleins\lh]=0$.
\end{rem}

Another important class of left-invariant complex structures are so-called abelian complex structures, which are characerised by $[\einsnull\lh,\einsnull\lh]=0$  or, equivalently, $[Jx,Jy]=[x,y]$  for all $x,y\in \lh$. In some sense this is the opposite condition to being a complex Lie-algebra and their deformations have been studied in \cite{mpps06, con-fin-poon06}. As  we pointed out in the introduction, deformations behave much more nicely in this case.

\section{Dolbeault cohomology}

In this section we will describe how the Dolbeault cohomology  of a nilmanifold with left-invariant complex structure $(M,J)$  is completely controlled by the Lie-algebra with complex structure $(\lh, J)$. This reduces many problems in the study of nilmanifolds to finite dimensional linear algebra. We will soon concentrate on the complex parallelisable case.

Let $(M,J)$ be a nilmanifold with left-invariant complex structure and $\lh$ be the Lie-algebra of the corresponding Lie-group.

We can identify elements in 
\[\Lambda ^{p,q}:=\Lambda^{p,q}(\lh^*,J)=\Lambda^p\einsnull{\lh^*}\tensor\Lambda^q\nulleins{\lh^*}\]
with left-invariant differential forms of type $(p,q)$ on $M$. The differential $d=\del+\delbar$ restricts to 
\[\Lambda^*\lh_\IC^*=\bigoplus\Lambda ^{p,q}\]
 and can in fact be defined in terms of the Lie bracket only: for $\alpha \in \lh^*$ and $x,y\in \lh$ considered as differential form and vectorfields we have
\begin{equation}\label{differential}
 d\alpha(x,y)=x(\alpha(y))-y(\alpha(x))-\alpha([x,y])=-\alpha([x,y])
\end{equation}
since all left-invariant functions are constant.

Let $H^k(\lh, \IC)$ be the $k$-th cohomology group of the complex 
\[\Lambda^*\lh^*_\IC: \quad 0\to \IC\overset{0}{\to} \lh_\IC^*\overset{d}{\to}\Lambda^2\lh_\IC^* \overset{d}{\to}\Lambda^3\lh_\IC^* \overset{d}{\to}\dots\]
and $H^{p,q}(\lh,J)$ be the $q$-th cohomology group of the complex
\[\Lambda^{p,*}:\quad 0\to \Lambda^{p,0}\overset{\delbar}{\to}\Lambda^{p,1}\overset{\delbar}{\to}\Lambda^{p,2}\overset{\delbar}{\to}\dots\]

In fact, the first complex calculates the usual Lie-algebra cohomology with values in the trivial module $\IC$ while the second calculates the cohomology of the Lie-algebra $\nulleins{\lh}$ with values in the module $\Lambda^{p,0}$ (see \cite{rollenske07b}).

\begin{theo}\label{cohom}
Let $M=\Gamma\backslash H$ be a real nilmanifold with Lie-algebra $\lh$. 
\begin{enumerate}
 \item The inclusion of $\Lambda^*\lh^*_\IC$ into the de Rham complex induces an isomorphism
\[H_{\mathrm{dR}}^*(M, \IC) \isom H^*(\lh, \IC)\]
in cohomology. (Nomizu, \cite{nomizu54}) 
\item  The inclusion of $\Lambda^{p,*}$ into the Dolbeault complex induces an inclusion
\begin{equation}\label{iota}
 \iota_J:H^{p,q}(\lh,J)\to H^{p,q}(M,J)
\end{equation}
which is an isomorphism if $(M,J)$ is complex parallelisable (Sakane, \cite{sakane76}) or if $J$ is abelian (Console and Fino, \cite{con-fin01}). Moreover, there exists a a dense open subset $U$ of the space of all left-invariant complex structures on $M$ such that $\iota$ is an isomorphism for all $J\in U$ (\cite{con-fin01})
\end{enumerate}
\end{theo}
Other work in this direction was done by Cordero, Fern\'andes, Gray and Ugarte \cite{cfgu00}. Conjecturally $\iota$ is an isomorphism for all left-invariant complex structures; in particular no counterexample is known.

For further reference we describe some cohomology groups in these terms.
\begin{lem}
Let $\lg$ be a complex Lie-algebra. Let us denote by $K^k:=\im(d:\Lambda^{k-1}\lh_\IC^*{\to}\Lambda^{k}\lh_\IC^*)$ the space of $k$-boundaries. Then
\begin{gather*}
 H^0(\lg, \IC)=\IC,\\
H^1(\lg,\IC)=\Ann(\kc_1\lg)=\Ann([\lg,\lg]),\\
K^2=\Ann(\ker([-,-]:\Lambda^2\lg\to \lg)).
\end{gather*}
Moreover, $H^{0,1}(\lg)=\overline{H^1(\lg,\IC)}$ and $\im(\delbar:\bar\lg^*\to \Lambda^2\bar\lg^*)=\bar K^2$.
\end{lem}
\pf
All assertions follow immediately from the fact that the differential $d:\lg^*\to\Lambda^2\lg^*$ is the dual of the Lie bracket $[-,-]:\Lambda^2\lg\to \lg$ and from the identification $\lg_\IC=\lg\oplus \bar \lg$.\qed

Since we are interested in deformations, the cohomology of the holomorphic tangent bundle (resp. tangent sheaf) $\Theta_{(M,J)}$ is of particular interest. It has been calculated in \cite{rollenske07b} for left-invariant complex structures for which \refb{iota} is an isomorphism, generalising results on abelian complex structure in \cite{mpps06, con-fin-poon06}.

But for a complex parallelisable nilmanifold $X$ we can calculate it directly (as observed by Nakamura \cite{nakamura75}).
Any element of the complex Lie-algebra $\lg$ gives rise to a holomorphic vector field. Hence the tangent sheaf is isomorphic to $\ko_X\tensor \lg$ and in cohomology we have a natural isomorphism
\[ H^q(X,\Theta_X)=H^q(X, \ko_X\tensor \lg)\isom  H^q(X, \ko_X)\tensor \lg =H^{0,q}(X)\tensor \lg\isom H^{0,q}(\lg)\tensor \lg.\]

Combining this with the previous results we get
\begin{lem}\label{cohomcalc}
Let $X=\Gamma\backslash G$ be a complex parallelisable nilmanifold. Then the tangent sheaf $\Theta_X\isom \ko_X\tensor \lg$ and its cohomology is calculated by the complex
\[ 0\to \lg \overset{0}{\to} \bar{\lg}^*\tensor \lg\overset{\delbar}{\to}\Lambda^{2} \bar{\lg}^*\tensor \lg\overset{\delbar}{\to}\dots\]
where the differential of $\bar\alpha\tensor X \in \Lambda^{p,0}\lg$ is given by $ \delbar (\bar\alpha\tensor X)=(\delbar\bar\alpha) \tensor X$.

In particular we have
\begin{gather*}
 H^0(X, \Theta)=\lg\\
H^1(X,\Theta)=H^1(X, \ko_X)\tensor \lg=\overline {\Ann([\lg,\lg])}\tensor \lg
\end{gather*}
\end{lem}

\section{Kuranishi theory}\label{kurani}

In \cite{kuranishi62} Kuranishi showed that for every compact complex manifold $X$ there exists a locally complete family of deformations which is versal at $X$. He constructs this family explicitly as a small neighbourhood of zero in the space of harmonic $(0,1)$-forms with values in the holomorphic tangent bundle after choosing some hermitian metric on $X$ (which always exists).

We will now apply his construction to complex parallelisable nilmanifolds using the results of the last section. 

Let $(M,J)=(\Gamma\backslash H,J)$ be the real nilmanifold with left-invariant complex structure underlying a complex parallelisable nilmanifold $X=\Gamma\backslash G$. The complex structure $J:\lh\to \lh$ is uniquely determined by the eigenspace decomposition $\lh_\IC=\einsnull\lh\oplus \nulleins\lh$.

A (sufficiently small) deformation of this decomposition $\lh_\IC=V\oplus \bar V$ can be encoded in a map $\Phi:\nulleins\lh \to \einsnull \lh$ such that $\bar V=(\id+\Phi) \nulleins \lh$, i.e., the graph of $\Phi$ in $\lh_\IC$ is the new space of vectors of type $(0,1)$.
This decomposition then determines a unique almost complex structure $J_V$ which is integrable if and only if $[V,V]\subset V$.

So far we have only described deformations of $J$ which remain left-invariant; this will be justified in a moment.

The integrability  condition is most conveniently expressed using the so-called Schouten bracket: for $X,Y\in \einsnull\lh$ and  $(0,1)$-forms $\bar \alpha, \bar\beta\in\nulleins{\lh^*}$ we set
\begin{equation}\label{Schouten}
 [\bar\alpha\tensor X, \bar\beta\tensor Y]:=\bar \beta \wedge L_{Y}\bar \alpha \tensor X+ \bar\alpha \wedge L_{X}\bar\beta\tensor Y+\bar \alpha\wedge \bar \beta \tensor [X,Y]
\end{equation}
 where $L_{X}\bar\beta=i_Xd\bar\beta+d(i_X \bar\beta)$ is the Lie derivative and $i_X$ is the contraction with $X$.

One can then show that the new complex structure is integrable if and only if $\Phi$ satisfies the Maurer-Cartan equation
\begin{equation}\label{MC}
 \delbar \Phi +[\Phi, \Phi]=0
\end{equation}
and it is well known that infinitesimal deformations, which correspond to  first-order solutions, are parametrised by classes in $H^1(X,\Theta_X)$ (see for example \cite{catanese88} or \cite{Huybrechts} for an overview).
But different solutions may well yield isomorphic deformations.

In order to single out a preferred solution we choose a hermitian structure on $\lg$ which induces a left-invariant hermitian structure on $X$.
Using the Hodge star operator associated to the hermitian metric we can define the formal adjoint $\delbar^*$ to $\delbar$ and the Laplace operator
\[\Delta:=\delbar\delbar^*+\delbar^*\delbar.\]

Defining the space of  harmonic forms to be $\kh^k=\ker (\Delta:\Lambda^{k} \bar{\lg}^*\to \Lambda^{k} \bar{\lg}^* )$ there is an orthogonal decomposition
\[\Lambda^{k} \bar{\lg}^*=B^k\oplus\kh^k\oplus V^k\]
where $B^k=\im(\delbar:\Lambda^{k-1} \bar{\lg}^*\to \Lambda^{k} \bar{\lg}^*)$ and $V^k=\im(\delbar^*:\Lambda^{k+1} \bar{\lg}^*\to \Lambda^{k} \bar{\lg}^*)$; this is just the intersection of  usual Hodge-decomposition with the subcomplex of left-invariant differential forms. The main point is that all harmonic forms are in left-invariant in our setting.

Since $\ker(\delbar)=B^k\oplus \kh^k$ we get an isomorphism
\[H^k(X,\Theta_X)\isom H^k(X,\ko_X)\tensor \lg \isom \kh^k\tensor \lg.\]

We are especially interested in the first two cohomology groups. By Lemma \ref{cohomcalc} we have $B^1=0$ which yields a commutative diagram
\[\xymatrix{
   \bar\lg^*\tensor \lg \ar[rr]^\delbar \ar@{=}[d]&& \ar@{=}[d]\Lambda^2\bar\lg^*\tensor \lg\\
(\kh^1\tensor\lg)\oplus (V^1 \tensor \lg)\ar[rr]^\delbar \ar[d]^{\mathrm{pr}} && (B^2\tensor \lg)\oplus( \kh^2\tensor \lg)\oplus( V^2\tensor \lg) \ar[dl]^P\ar[d]_H\\
V^1\tensor \lg & \ar[l]^{\delta}_{\isom} B^2\tensor\lg &\kh^2\tensor \lg.
}\]
We denote by $\delta$ the inverse of the isomorphism $P\circ \delbar: V^1\tensor\lg\to B^2\tensor\lg$.

We will now use these operators to describe the Kuranishi space: let $X_1,\dots, X_n$ be a basis of $\lg$ and $\bar\omega^1, \dots \bar\omega^m$ be a basis for $\kh^1$. Then $\{\bar\omega^i\tensor X_j\}$ is a basis of $H^1(X,\Theta_X)$ and we define recursively
\begin{equation}\label{Phi}
 \begin{split}
 \Phi_1(\underline t)&=\sum_{i,j} t_i^j \bar\omega^i\tensor X_j, \\
\Phi_2(\underline t)&:=-
\delta\circ P [\Phi_1(\underline t), \Phi_1(\underline t)],\\
\Phi_k(\underline t)&:=-
\delta\circ P \sum_{1\leq i<k} \left[\Phi_i(\underline t), \Phi_{k-i}(\underline t)\right]  \quad (k\geq 2),
 \end{split}
\end{equation}
obtaining a formal power series
\[\Phi(\underline t)=\sum_{k\geq1} \Phi_k(\underline t).\]
 \margincom{show or cite that $\sum_{1\leq i<k} [\Phi_i, \Phi_{k-i}] $ is closed!(graded Lie-algebra structure)}

We see that $\Phi_k$ is a homogeneous polynomial  of degree $k$ in the variables $t_i^j$ and it is easy to verify that 
\[\delbar\Phi+[\Phi,\Phi]=H[\Phi,\Phi].\]
The map $\Phi$ does not depend on the choice of the basis and we can define the obstruction map
\[ \obs:\kh^1\tensor \lg \to \kh^2\tensor \lg,\qquad \mu=\sum_{i,j} t_i^j \bar\omega^i\tensor X_j\mapsto H[\Phi(\underline t),\Phi(\underline t)].\]

We can now formulate Kuranishi's theorem in our context.

\begin{theo}[\cite{kuranishi62}] 
The formal powerseries $\Phi(\underline t)$ converges for sufficiently small values of $\underline t$ and there is a versal family of deformations of $X$ over the space
\[\Kur(X):=\{\mu\in \kh^1(\Theta_X)\mid  \|\mu\|<\epsilon; \obs(\mu)=0\}.\]
where $\kh^1(\Theta_X)=\kh^1\tensor \lg$ is the space of harmonic 1-forms with values in $\Theta_X$. $\Kur(X)$ is called the Kuranishi space of $X$.
\end{theo}

By construction $\Phi$ is left-invariant and hence the new complex structure will also be left-invariant. In fact, the new subbundle of tangent vectors of type $(0,1)$ in $TM_\IC$ is obtained by translating the subspace $(\id+\Phi)\nulleins\lg\subset \lg_\IC$.
We have reproved that all sufficiently small deformations of our complex parallelisable nilmanifold carry a left-invariant complex structure.

Note that the construction involved the choice of a hermitian structure so $\Kur(X)$ is not defined in a canonical way. Nevertheless for different choices of a metric (the germs of) the resulting spaces are (non canonically) isomorphic.

The values of $\underline t$ have to be small for two different reasons. First of all we need to ensure the convergence of the formal power series $\Phi(\underline t)$ and secondly $(\id+\Phi)\bar \lg$ should be the space of $(0,1)$ vectors for an integrable almost complex structure, in other words we need $(\id+\Phi)\bar \lg\oplus\overline{(\id+\Phi)\bar \lg}=\lg_\IC$. We will see that the first issue will not arise in our setting.

\begin{rem}
Usually the terms of the formal power series $\Phi$ are described using Green's operator, which inverts the Laplacian on the orthogonal complement of harmonic forms, setting
\[\Phi_k(\underline t):=-\delbar^*G \sum_{1\leq i<k} \left[\Phi_i(\underline t), \Phi_{k-i}(\underline t)\right]\]
It is straight-forward to check that this agrees with our definition above using the identities $G\circ\Delta+H=\Delta\circ G+H=\id$ and definition of the Laplacian.

Our formula involves only $\delta=\inverse \delbar$ and the projection $P$ which will facilitate the computation of examples in Section \ref{examples}.
\end{rem}

Now that we have seen how the Kuranishi space is constructed we want to investigate its structure in detail for complex parallelisable nilmanifolds.

The key result is the following:

\begin{lem}\label{schoutennil}
 Let $\bar\alpha\tensor X, \bar\beta\tensor Y \in \bar\lg^*\tensor \lg$. Then their Schouten bracket is
\[[\bar\alpha\tensor X, \bar\beta\tensor Y ]=\bar\alpha\wedge \bar \beta\tensor[X,Y].\]
\end{lem}
\pf Comparing the expression with the general formula \refb{Schouten}  it suffices to show that for $X\in \lg$ and $\bar\alpha\in\bar\lg$ the Lie-derivative $L_X\bar\alpha=i_Xd\bar\alpha+d(i_X \bar\alpha)=0$. But $\bar\alpha$ is of type $(0,1)$ and $d\bar\alpha$ is of type $(0,2)$ (since $[\lg, \bar\lg]=0$) so both vanish when contracted with a vector of type $(1,0)$. \qed

This gives us 

\begin{lem} For $\Phi$ as in the recursive description \refb{Phi} we have
 $[\Phi_k, \Phi_l] \in \Lambda^2\bar{\lg}^*\tensor \kc_{k+l-1}\lg\subset \Lambda^2\bar{\lg}^*\tensor\lg$.
\end{lem}
\pf We prove our claim by induction: for $k=1$ there is nothing to prove since $\kc_0\lg=\lg$. Note that, by the Jacobi identity, $[\kc_k\lg, \kc_l\lg]\subset \kc_{k+l+1}\lg$. Since the Schouten bracket is  the Lie bracket on the vector part and the map $\delta=\inverse\delbar$ acts only on the form part our claim follows. \qed

We deduce immediately that the Kuranishi space can not be too complicated:
\begin{theo}\label{polynom}
If $\lg$ is $\nu$-step nilpotent and $\Phi$ as in \refb{Phi} then 
\[\obs(\underline t)=\sum_{\stackrel{1\leq i,j< \nu,}{ i+j\leq \nu}}H[\Phi_i, \Phi_j].\]
In particular $\Kur(X)$ is cut out by polynomial equations of degree at most $\nu$.
\end{theo}

\pf Since $\lg$ is $\nu$-step nilpotent $\kc_k\lg=0$ for $k\geq \nu$. By the previous Lemma this implies that $[\Phi_i, \Phi_j]=0$ whenever $i+j-1\geq \nu$ and hence  the only possibly non-vanishing terms of $\obs=H[\Phi,\Phi]$ are the ones given above.\qed

For further reference we use Lemma \ref{schoutennil} to calculate the second order obstructions, i.e., the quadratic term of the obstruction map $\obs$: let as before $\bar\omega^1, \dots, \bar\omega^{m}$ be a basis of $\kh^1=\overline {\Ann([\lg,\lg])}$ and $X_1, \dots, X_n$ be a basis of $\lg$.

Then we can represent any element in $H^1(X, \Theta_X)$ as 
\[ \Phi_1(\underline t)=\sum_{i,j} t_i^j \bar\omega^i\tensor X_j\]
and consequently
\begin{equation}\label{deg2}
 \begin{split}
[\Phi_1(\underline t),\Phi_1(\underline t)]&=[\sum_{i,k} t_i^k \bar\omega^i\tensor X_k, \sum_{j,l} t_j^l \bar\omega^j\tensor X_l]\\
  &= \sum_{i,j,k,l}(t_i^k t^l_j)[ \bar\omega^i\tensor X_k, \bar\omega^j\tensor X_l]\\
  &= \sum_{i,j,k,l}(t_i^k t^l_j)\bar\omega^i\wedge\bar\omega^j\tensor [X_k, X_l]\\
& = \sum_{1\leq i<j\leq m}\sum_{k,l}(t_i^k t^l_j-t_j^k t^l_i)\bar\omega^i\wedge\bar\omega^j\tensor [X_k, X_l]\\
& = \sum_{1\leq i<j\leq m}\sum_{1\leq k<l\leq n}2(t_i^k t^l_j-t_j^k t^l_i)\bar\omega^i\wedge\bar\omega^j\tensor [X_k, X_l]\\
& = 2 \sum_{1\leq i<j\leq m}\sum_{1\leq k<l\leq n}
\det\begin{pmatrix}t_i^k & t_i^l\\t_j^k & t_j^l \end{pmatrix}
\bar\omega^i\wedge\bar\omega^j\tensor [X_k, X_l].
 \end{split}
\end{equation}

We deduce from this formula a necessary condition for the Kuranishi space to be smooth:
\begin{lem}\label{lambda2}
If  the subspace  $\Lambda^2 \kh^1 \subset \Lambda^2 \bar\lg$ is not contained in $B^2$  and $\lg$ is not abelian then there is a non-vanishing obstructions in degree 2 and the Kuranishi space is singular.
\end{lem}
\pf
Assume that $\Lambda^2 \kh^1 \subset \Lambda^2 \bar\lg$ is not contained in $B_2$. Then there is some basis vector $\bar\omega^i\wedge\bar\omega^j$ which is not contained in the image of $\delbar$. Since $\lg$ is not abelian there are vectors $X_k, X_l$ such that  $[X_k, X_l]\neq 0$. Setting $t_p^q=0$ if $p\neq i,j$ or $q\neq k,l$ and choosing the remaining coefficient such that $\det\begin{pmatrix}t_i^k & t_i^l\\t_j^k & t_j^l \end{pmatrix}\neq 0$ we have found an obstructed element in $H^1(X, \Theta_X)$. \qed

The condition that the Kuranishi space be smooth is very strong. To make this more precise we need to recall the definition of  the free 2-step Lie-algebra: let $m\geq2$, $V=\IC^m$ and $\gothb_m:=V\oplus \Lambda^2 V$. Then $\gothb_m$ with the Lie bracket
\[ [ \cdot, \cdot]: \gothb_m\times \gothb_m\to \gothb_m, \qquad [a+b\wedge c, a'+b'\wedge c']:=a\wedge a'\]
 is the free 2-step nilpotent Lie-algebra.

\begin{theo}\label{freecond}
If $\lg$ is not abelian then there is a non-vanishing obstruction in degree 2 if and only if $\lg\slash \kc^2\lg$ is not isomorphic to a free 2-step nilpotent Lie-algebra.

Hence, if $\lg\slash \kc^2\lg$ is not free then the Kuranishi space is singular.
\end{theo}

The vanishing of all obstructions on degree 2 is not a sufficient condition for the Kuranishi space to be smooth. A 4-dimensional example where $\Kur(X)$ is cut out by a single cubic equation can be found in Section \ref{explicit4}.

All examples with smooth Kuranishi space which we could find were actually free Lie algebras and at least in the 2-step nilpotent case there are no other:

\begin{cor}\label{b_m}
If $\lg$ is 2-step nilpotent then the Kuranishi space is smooth if and only if $\lg$ is a free 2-step nilpotent Lie-algebra, i.e., $\lg\isom \gothb_m$ with $m=h^{0,1}(X)$.
\end{cor}
\pf This follows immediately from the theorem since for a 2-step nilpotent Lie-algebra we have $\kc_2\lg=0$, hence $\lg/\kc_2\lg\isom\lg.$\qed

Note that $\gothb_2$ is the complex Heisenberg algebra, which is the Lie-algebra of the universal cover of the Iwasawa-manifold. So we have reproved the smoothness of the  Kuranishi space of the Iwasawa manifold first observed by Nakamura.	

It is very easy to produce examples with singular Kuranishi space:

\begin{cor}
If  $\lg\isom \lg'\oplus \gotha$ where $\gotha\isom \IC^n$ is an abelian Lie-algebra and $\lg'$ is not abelian, then the Kuranishi-space is singular.

In particular, if $X$ is any complex parallelisable nilmanifold which is not a torus and $T$ is a complex torus then $X\times T$ has obstructed deformations.
\end{cor}
\pf We have $\lg\slash \kc_2\lg=\lg'\slash\kc_2\lg'\oplus \gotha$ which is not free. An application of the theorem proves the assertion.\qed

Before we can address the proof of Theorem \ref{freecond} we need a technical lemma.

Let $\lg=\kc_0 \lg\supset \kc_1\lg\supset \dots \supset \kc_{\nu}\lg=0$ be the descending central series and let
$\kc^k\lg^*=\Ann\kc_k\lg$. We get a filtration
\[0=\kc^0\lg^*\subset \kc^1\lg^*=\Ann(\kc_1\lg)\subset \dots\subset \kc^\nu\lg^*=\lg^*.\]

\begin{lem}\label{dck}
Setting
\[ W^k=\langle \alpha\wedge \beta \in \Lambda^2\lg^*\mid \alpha \in \kc^i\lg^*, \beta\in \kc^j\lg^*, i+j\leq k\rangle 
\subset \Lambda^2\lg^*\]
we have
\[ d\alpha \in W^k \iff \alpha \in \kc^k\lg^*.\]
\end{lem}
\pf
 Assume that there is $\alpha \notin \kc^k\lg^*$ with $d\alpha \in W^k$.
By the Jacobi-identity $\kc_k\lg$ is generated by elements of the form $X=[Y,Z]$ where $Y\notin \kc_1\lg$ and $Z\in \kc_{k-1}\lg$ and hence $\alpha(X)\neq 0$ for one such element.
By the definition of $\kc^i\lg^*$ we have $\beta(Y,Z)=0$ for all $\beta \in W^k$. On the other hand 
\[d\alpha(Y,Z)=-\alpha([Y,Z])=-\alpha(X)\neq 0\] so  $d\alpha\notin W_k$ -- a contradiction.

The other direction is a well known fact for nilpotent Lie algebras. It can be easily seen picking a  basis adapted to the descending central series (often called Malcev or Engel basis) and writing $\alpha$ as a linear combination of the elements of  the dual basis.
\qed

\emph{Proof of Theorem \ref{freecond}.} 
Let $\lg$ be a non-abelian Lie-algebra. By Lemma \ref{lambda2} it suffices to show that $\Lambda^2\kh^1\subset B_2$ if and only if $\lg\slash \kc^2\lg$ is a free 2-step Lie-algebra.

Recalling that $\kh^1=\overline{\kc^1\lg^*}$ and $B^2=\overline{\im(d)}$  we have to prove that $\Lambda^2\kc^1\lg^*=W^1$ is in the image of the differential if and only if $\lg\slash \kc^2\lg$ is free.

The Lie bracket in $\lg$ can also be considered as a linear map
\[b:\Lambda^2 \lg \to \kc_1\lg,\]
which is, by definition, surjective.
Dualising we get (the restriction of) the differential
\[  d: (\kc_1\lg)^* \to \Lambda^2 \lg^*,\]
which is now injective.

Let $A$ be the anullator of $\kc_2\lg$ in $(\kc_1\lg)^*$. Then we infer from Lemma \ref{dck}
that $d\restr{A}:A\to W^2=\Lambda^2\kc^1\lg^*$, in fact, 
\[ dA=\im (d) \cap W^2=\im(d)\cap \Lambda^2\kc^1\lg^*.\]
The dual map
\[b': \left(\Lambda^2\kc^1\lg^*\right)^*=\Lambda^2(\lg\slash\kc_1\lg) \to A^*=\kc_1\lg\slash\kc_2\lg\]
gives an anti-symmetric bilinear form on $\lg\slash\kc_1\lg$ with values in $\kc_1\lg\slash\kc_2\lg$ which is exactly the Lie bracket in the quotient Lie-algebra $\lg\slash\kc_2\lg$.

Hence we see that $\Lambda^2\kc^1\lg^*$ is in the image of $d$ if and only if 
$d: A \to W^2$
is surjective if and only if $b'$ is injective. But $b'$ is by definition surjective so it is injective if and only if it is bijective in which case $\Lambda^2(\lg\slash\kc_1\lg)\isom\kc_1\lg\slash\kc_2\lg$
and the Lie-algebra $\lg\slash \kc^2\lg$ is indeed free. \qed

\section{Deformations remaining complex parallelisable}\label{remainparall}

It is a natural question if there are conditions which guarantee that a given small deformation of our complex parallelisable manifold $X$ is again complex parallelisable.

So let $\mu \in H^1(X, \Theta_X)=\kh^1\tensor \lg$ be a infinitesimal deformation and $\Phi$ the corresponding iterative solution of the Maurer-Cartan equation as in \refb{Phi}. The new space of $(0,1)$-vectors is $(\id+\Phi)\bar \lg$. (Recall that we identified $\lg_\IC=\lg\tensor\bar\lg$.)

By Lemma \ref{parallchar}  the new complex structure is again complex parallelisable if and only if
\[ [(\id+\bar\Phi) X, (\id+\Phi)\bar Y]=0\]
for all $X,Y\in \lg$.
Looking at the terms up to first order yields
\[ [X, \bar Y]+[\bar\mu X,\bar Y] +[ X,\mu\bar Y]=[\bar\mu X,\bar Y] +[ X,\mu\bar Y]=0.\]
The first of these terms is in $\bar\lg$ while the second is in $\lg$ and they are complex conjugate to each other up to sign and renaming. Thus we call $\mu$ an infinitesimally complex parallelisable deformation if
\begin{gather*}
\forall X,Y \in \lg : [ X,\mu\bar Y]=0 
\iff \mu \in \kh^1\tensor \kz\lg.
\end{gather*}
 Such infinitesimal deformations are always unobstructed:
if $\mu \in \kh^1\tensor \kz\lg$ then $[\mu,\mu]\in \Lambda^2\bar\lg^*\tensor [\kz\lg, \kz\lg]=0$. Hence in the recursive definition \refb{Phi} all higher order terms vanish, $\Phi=\mu$ and $\obs(\mu)=0$.

We have proved
\begin{theo}\label{remaining}
For an element $\mu \in H^{1}(X, \Theta_X)=H^1(X,\ko_X)\tensor \lg$ the following are equivalent:
\begin{enumerate}
\item $\mu \in H^1(X, \ko_X)\tensor \kz\lg$.
 \item $\mu$ defines an infinitesimally complex parallelisable deformation.
\item $t\mu$ induces a 1-parameter family of complex parallelisable manifolds for $t$ small enough, i.e., provided that $(\id+t\mu)\bar\lg\oplus(\id+t\bar\mu)\lg=\lg_\IC$.
\end{enumerate}
Hence the Kuranishi family is (locally) a cylinder over an analytic subset of $H^1(X, \ko_X)\tensor (\lg/\kz\lg)$.
\end{theo}

Since $\lg=\kz\lg$ if and only if $\lg$ is abelian we deduce:
\begin{cor}\label{nondef}
 If $\lg$ is not abelian then there are small deformations of $X$ which are not complex parallelisable.
\end{cor}

\section{Examples}\label{examples}

 
We continue to use the notation from Remark \ref{notation}. The deformation theory of the complex parallelisable nilmanifold $X$ is completely determined by the Lie-algebra $\lg$ and we have already discussed two series of examples where $\Kur(X)$ is smooth.
\begin{itemize}
\item If $\lg=\mathfrak a_k$ is the $k$-dimensional abelian Lie-algebra then $X$ is a torus and $\Kur(X)$ is smooth of dimension $\frac{k^2(k+1)}{2}$.
\item If $\lg=\mathfrak b_m$ is the free 2-step nilpotent Lie-algebra on $m$ generators, which has dimension  $\frac{m(m+3)}{2}$,  then $\Kur(X)$ is smooth of dimension $\frac{m^2(m+3)}{2}$ (see Corollary \ref{b_m}).
\end{itemize}

\subsection{Examples in low dimension -- overview}

Nilpotent complex Lie algebras are classified up to dimension 7 \cite{magnin86} and partial results are known in dimension 8 \margincom{add reference for dim8}. Starting from dimension 7 there are infinitely many non-isomorphic cases.

We will now describe the Kuranishi-space of complex parallelisable nilmanifolds up to dimension 5. 

There is a convenient way to describe a nilpotent Lie-algebra $\lg$ using the differential $d: \lg\to \Lambda^2\lg$. The expression
\[ \lg=(0,0,0,0,12+34)\]
means the following: with respect to a basis $\omega^1, \dots , \omega^5$ the differential is given by
\[ d\omega^1=d\omega^2=d\omega^3=d\omega^4=0 \text{ and } d\omega^5=\omega^1\wedge\omega^2+\omega^3\wedge\omega^4.\]
This determines the Lie bracket, which is the dual map (see \refb{differential}).

More precisely, if we denote by  $X_1, \dots, X_5$ the dual basis then the only non-zero Lie brackets are $[X_1, X_2]=[X_3,X_4]=-X_5$.

Table 1 lists all Lie-algebras up to dimension 5 in this notation together with some  information on the Kuranishi space of an associated complex parallelisable nilmanifold. We denote the nilpotency index by $\nu$.

Note that all Lie-algebras with smooth Kuranishi space are either free or abelian. One can check that also the free 4-step nilpotent Lie-algebra on 2 generators $(0,0,12,13,23,14,25,24+15)$ has smooth Kuranishi space.
 \begin{table}\label{alle}
 \caption{Kuranishi spaces up to dimension 5.}
\begin{center}
\begin{tabular}{ccccccc}
 $\dim$ & Lie-algebra & $\nu$ &$h^1(\Theta_X)$ & smooth & irreducible & reduced\\ \hline
1 & $\gotha_1$ &1& 1 & $\surd$ & $\surd$ & $\surd$ \\ \hline
2 & $\gotha_2$ &1& 6 & $\surd$ & $\surd$ & $\surd$ \\ \hline
3 & $\gotha_3$ &1& 18 & $\surd$ & $\surd$ & $\surd$ \\ 
3 & $\gothb_1$ &2& 6 & $\surd$ & $\surd$ & $\surd$ \\ \hline
4 & $\gotha_4$ &1& 40 & $\surd$ & $\surd$ & $\surd$ \\
4 & $(0,0,0,12)$ &2& 12 & $-$ & $-$ & $\surd$ \\ 
4 & $(0,0,12,13)$ &3& 8 & $-$ & $-$ & $\surd$ \\ \hline
5 & $\gotha_5$ &1& 75 & $\surd$ & $\surd$ & $\surd$ \\ 
5 & $(0,0,0,12,13)$ &2& 15 & $-$ & $-$ & $\surd$ \\ 
5 & $(0,0,0,0,12+34)$ &2& 20 & $-$ & $-$ &  $\surd$ \\ 
5 & $(0,0,12,13,23)$ &3& 10 & $\surd$ & $\surd$ & $\surd$ \\
5 & $(0,0,0,12,13+24)$ &3& 15 & $-$  & $-$  & $-$ \\ 
5 & $(0,0,12,13,14)$ &4& 10 & $-$ & $-$ & $-$\\ 
5 & $(0,0,12,13,14+23)$ &4& 10 & $-$ & $-$  & $-$  \\
\end{tabular}
\end{center}

\leftline{\scriptsize ($\nu$ = nilpotency index)}
 \end{table}

\subsection{Examples in low dimension -- explicit descriptions}

In this section we will give explicit equations for the Kuranishi space of some examples.
In order to avoid cumbersome notation we will only consider the germ of the Kuranishi space at zero which will be denoted by $\Kur(X)_0$.

Since nothing interesting happens in dimension 1, 2, and 3 we start in dimension 4.

\subsubsection{Computations in dimension 4}\label{explicit4}

We will now compute the Kuranishi space explicitly  for the two singular examples in dimension 4. 

The structure equations of the considered Lie-algebras are given with respect to the bases $X_1, \dots, X_n$ and $\omega^1, \dots, \omega^n$ as described at the beginning of this section.  Thus we will always start the computation of the iterative solution of the Maurer-Cartan equation with the element
\[ \Phi_1(\underline t)=\sum_{i=1}^{n}\sum_{j=1}^{m} t_i^j \bar\omega^i\tensor X_j\]
where $n=\dim\lg$ and ${m}=\codim\kc_1\lg=h^{0,1}(X)$.

In order to use harmonic forms we equip $\lg$ with the unique hermitian metric such that the $X_i$ form an orthonormal basis.

In every step of the recursion \refb{Phi} we will decompose $[\Phi_k, \Phi_l]=\beta+\chi$ where $\chi$ is harmonic and $\beta$ is exact. Then $\chi$ will contribute to the obstruction map and $\delta(\beta)=\inverse{(\delbar)}\beta$ will, if necessary, be used to compute the next iterative step.

\paragraph*{The Lie-algebra $\lg=(0,0,0,12)$}

Since $\lg$ is 2-step nilpotent we only have to look at obstructions in degree 2, i.e., $\obs=H[\Phi_1, \Phi_1]$. Since $[X_1,X_2]=-X_4$ is the only non-zero bracket we deduce from \refb{deg2} that 
\begin{align*}
[\Phi_1(\underline t),\Phi_1(\underline t)]&= - 2 \sum_{1\leq i<j\leq 3}\det\begin{pmatrix}t_i^1 & t_i^2\\t_j^1 & t_j^2 \end{pmatrix}
\bar\omega^i\wedge\bar\omega^j\tensor X_4\\
&=-2\det\begin{pmatrix}t_1^1 & t_1^2\\t_3^1 & t_3^2 \end{pmatrix}
\bar\omega^1\wedge\bar\omega^3\tensor X_4-2\det\begin{pmatrix}t_2^1 & t_2^2\\t_3^1 & t_3^2 \end{pmatrix}
\bar\omega^2\wedge\bar\omega^3\tensor X_4\\
&\qquad-\delbar\left(2\det\begin{pmatrix}t_1^1 & t_1^2\\t_2^1 & t_2^2 \end{pmatrix}
\bar\omega^4\tensor X_4\right).
\end{align*}
Hence
\begin{align*}
\Kur(X)_0&=\{ \underline t \in \IC^{12}\mid \det\begin{pmatrix}t_1^1 & t_1^2\\t_3^1 & t_3^2 \end{pmatrix}=\det\begin{pmatrix}t_2^1 & t_2^2\\t_3^1 & t_3^2 \end{pmatrix}=0\}_0\\
&=\left(\IC^6\times Y\right)_0
\end{align*}
where 
\[Y= \{ t_3^1 = t_3^2=0\}\cup \{ \rk\begin{pmatrix}t_1^1 &t_2^1 & t_3^1 t_1^2\\t_1^2 & t_2^2& t_3^2 \end{pmatrix}\leq 1\}.\]
In particular we see that the Kuranishi space is a cylinder over the reducible space $Y$.

\paragraph*{The Lie-algebra $\lg=(0,0,12,13)$}

We infer from \refb{deg2} that
\begin{align*}
[\Phi_1(\underline t),\Phi_1(\underline t)]&=
- 2 \det\begin{pmatrix}t_1^1 & t_1^2\\t_2^1 & t_2^2 \end{pmatrix} \bar\omega^1\wedge\bar\omega^2\tensor X_3
- 2 \det\begin{pmatrix}t_1^1 & t_1^3\\t_2^1 & t_2^3 \end{pmatrix} \bar\omega^1\wedge\bar\omega^2\tensor X_4\\
&=-\delbar\left(2\det\begin{pmatrix}t_1^1 & t_1^2\\t_2^1 & t_2^2 \end{pmatrix} \bar\omega^3\tensor X_3
+ 2 \det\begin{pmatrix}t_1^1 & t_1^3\\t_2^1 & t_2^3 \end{pmatrix}\bar\omega^3\tensor X_4\right)
\end{align*}
 and by the recursion formula we set
\[\Phi_2:=2\det\begin{pmatrix}t_1^1 & t_1^2\\t_2^1 & t_2^2 \end{pmatrix} \bar\omega^3\tensor X_3
+ 2 \det\begin{pmatrix}t_1^1 & t_1^3\\t_2^1 & t_2^3 \end{pmatrix}\bar\omega^3\tensor X_4.\]

We see that there are no obstructions of second order and calculate (noting that $X_4$ is in the centre and that $[X_2, X_3]=0$)
\begin{align*}
[\Phi_1(\underline t),\Phi_2(\underline t)]&=[t^1_1\bar\omega^1\tensor X_1 + t^1_2 \bar\omega^2\tensor X_1, 2\det\begin{pmatrix}t_1^1 & t_1^2\\t_2^1 & t_2^2 \end{pmatrix} \bar\omega^3\tensor X_3]\\
&= -2 \det\begin{pmatrix}t_1^1 & t_1^2\\t_2^1 & t_2^2 \end{pmatrix} \left( t^1_1\bar\omega^1\wedge \bar\omega^3\tensor X_4 +t^1_2\bar\omega^2\wedge \bar\omega^3\tensor X_4\right)\\
&= -2\det\begin{pmatrix}t_1^1 & t_1^2\\t_2^1 & t_2^2 \end{pmatrix} \left( t^1_2\bar\omega^2\wedge \bar\omega^3\tensor X_4+t^1_1\delbar \bar\omega^4\tensor X_4 \right)\\
&=-2 t^1_2\det\begin{pmatrix}t_1^1 & t_1^2\\t_2^1 & t_2^2 \end{pmatrix}  \bar\omega^2\wedge \bar\omega^3\tensor X_4 \mod B_2
\end{align*}
Hence we have
\[\Kur(X)_0=\{ \underline t \in \IC^2\tensor \IC^4 =\IC^8\mid   t^1_2\det\begin{pmatrix}t_1^1 & t_1^2\\t_2^1 & t_2^2 \end{pmatrix}=0\}_0,\]
in other words, $\Kur(X)_0$ is a cylinder
 over the cone over the union of a plane and a quadric in $\IP^3$.

\subsubsection{Remarks on dimension 5}\label{dim5sect}

The computations in dimension 5 proceed along the same lines as in dimension 4 but are, as one might imagine,  much more involved. Thus, we will only present the results.

In the view of Theorem \ref{remaining} the Kuranishi space is a cylinder over an analytic subset of the vector space $H^1(\bar \lg,\IC)\tensor (\lg\slash \kz\lg)$ whose dimension we denote by $d$.

The germ of the Kuranishi space at 0 is cut out by polynomial function and we will give the primary decomposition, computed using the program Singular \cite{GPS05},  of the ideal $I$ of all these functions. Different ideals in the decomposition correspond to different irreducible components. 

If some component is set-theoretically contained in another component we call it an embedded component; this can only happen if the component is not reduced. Non-reduced components occur if there are  infinitesimal deformations which can be lifted up to a certain order but not to actual deformations.

In all examples the Kuranishi space has several irreducible components. We denote by $k$ be the number of components of the reduced space and by $e$ the number of embedded components. 
Note that in the case $\lg=(0,0,12,13,14+23)$ there are two non-reduced components which are not embedded, both supported on linear subspaces.

To simplify the description of the ideals we introduce the notation
\begin{gather*}
 \delta_{ij}^{kl}:=\det\begin{pmatrix}t_i^k & t_i^l\\t_j^k & t_j^l \end{pmatrix},\\
\Delta_{ijk}^{lmn}:=\det\begin{pmatrix}
t_i^l &t_j^l &t_k^l\\
t_i^m &t_j^m &t_k^m\\
t_i^n &t_j^n &t_k^n
\end{pmatrix}.
\end{gather*}
The results can now be found in Table 2 where we also give the codimension and the degree of the various components.

\begin{table}\label{dim5}
 \caption{Singular deformation spaces in dimension 5 (See Section \ref{dim5sect} for notations.)}
\begin{sideways}
 \begin{tabularx}{\textheight}{ccccccX}
$\mathbf{\lg}$ & $\mathbf d$ & $\mathbf{(k,e)}$ & \textbf{Codim.} &\textbf{Degree} & \textbf{reduced?} & \textbf{Ideal (primary decomposition)}\\
\hline &\qquad\\
$(0,0,0,12,13)$ &9&$(2,0)$&$(2,2)$&$(3,1)$& $(\surd,\surd)$&$(\delta_{23}^{23}, \delta_{23}^{13}, \delta_{23}^{12})\cap (t_3^1, t_2^1)$ \\
&\qquad\\
$(0,0,0,12,13+24)$ & 12 &$(3,2)$ & $(4, 5, 4);( 5, 5)$ &$(9, 3, 3);(2, 4)$ & $(\surd, \surd, \surd)$&
\begin{minipage}{6cm}
{
\begin{align*}
 &(\delta^{13}_{23}+\delta^{24}_{23}, \delta_{23}^{12}, \delta_{13}^{14}+\delta_{13}^{24}, \delta_{13}^{12}, \delta_{12}^{13}+\delta_{12}^{24}, \delta^{12}_{12})\\ 
 \cap &( t_3^2, t_3^1,t_1^2, t_2^1t_3^3+t^2_2t_3^4, 2(t^2_2)^2-t_3^3, 2t_2^1t_2^2+t_3^4)\\
\cap&( t_3^2, t_3^1, t_2^1t_3^3+t^2_2t_3^4, t^1_1 + t_1^2t_3^4, \delta_{12}^{12})\\
\cap&(t_3^2, t_1^2, (t_3^3)^2, t_3^1t_3^3, t_2^2t_3^3, (t_3^1)^2, \delta_{23}^{13}+t_2^2t_3^4, t_2^2t_3^1, \delta_{13}^{13}, (t_2^2)^2)\\
\cap&(t_1^2, t_3^1t_3^3+t_3^2t_3^4,  (t_3^2)^2, t_3^1t_3^2, t_1^1t_3^2,(t_3^1)^2, \\
&\qquad\delta^{13}_{23}+\delta^{24}_{23}, \delta_{23}^{12}, 
t_1^1t_3^1,(t_1^1)^2, \delta^{13}_{13}+t_1^4t_3^2+2t_1^1(t_2^2)^2)
\end{align*}}
\end{minipage}\\
&\qquad\\
$(0,0,12,13,14)$ & 8 & $(2,1)$ & $(2,1);(2)$ & $(3,1);(2)$ & $(\surd,\surd)$ &
\begin{minipage}{4cm}
\[(\delta_{12}^{23}, \delta_{12}^{13}, \delta_{12}^{12})\cap (t_2^1)\cap (\delta^{12}_{12}, t_1^1t_2^1, (t_1^1)^2, (t_2^1)^2)\]
\end{minipage}
\\
&\qquad\\
$(0,0,12,13,14+23)$ & 8 & $(4,0)$ & $(2, 2, 2, 2)$ & $(3,2,2,3)$ & $(\surd, \surd,-, -)$ & 
\begin{minipage}{8cm}
{\begin{align*}
 & (\delta_{12}^{23}, \delta_{12}^{13}, \delta^{12}_{12}) \cap (t_2^1, 2(t_1^1)^2-t_2^2)\\
\cap &((t_2^2)^2,t_2^1t_2^2, (t_2^1)^2, 2t_1^1\delta_{12}^{12}-t_2^1t_2^3)\\
 \cap& ( (t_2^1)^3, t_2^2\delta_{12}^{12}+t_2^1\delta_{12}^{13}, t_2^1\delta_{12}^{12},\\
& \qquad t_1^1(t_2^1)^2,t_1^2\delta_{12}^{12} +t_1^1\delta_{12}^{13}, t_1^1\delta_{12}^{12}, (t_1^1)^2t_2^1, (t^1_1)^3)
\end{align*}}
\end{minipage}\\
&\qquad\\
$(0,0,0,0,12+34)$ & 16 & $(2,0)$ & $(5,5)$ & $(20,12)$ & $(\surd, \surd)$ &
\begin{minipage}{8cm}
{\begin{align*}
(&\delta^{12}_{34}+\delta_{34}^{34}, \delta^{12}_{24}+\delta^{34}_{24}, \delta^{12}_{14}+\delta^{34}_{14}, \delta^{12}_{23}+\delta^{34}_{23}, \delta^{12}_{13}+\delta^{34}_{14},\\
& \delta^{12}_{12}+\delta^{34}_{12}, \Delta^{234}_{234}, \Delta^{234}_{134}, \Delta^{234}_{124}, \Delta^{134}_{234}, \Delta^{134}_{134}, \Delta^{134}_{124}, \Delta_{123}^{234}, \Delta^{134}_{123})\\
\cap&(\delta^{12}_{24}+\delta^{34}_{23}, \delta^{12}_{14}+\delta^{34}_{14}, \delta^{12}_{23}+\delta^{34}_{23}, \delta^{12}_{13}+\delta^{34}_{13}, \delta^{12}_{34}-\delta^{34}_{12},\\
 &\delta^{24}_{12}+\delta^{24}_{34}, \delta^{23}_{12}+\delta^{23}_{34}, \delta^{14}_{12}+\delta^{14}_{34}, \delta^{13}_{12}+\delta^{13}_{34},\delta^{12}_{12}-\delta^{34}_{34},\\
& t^1_2\delta^{12}_{34}-t^1_4\delta^{34}_{23}+t^1_3\delta^{14}_{34}, t^1_1\delta^{12}_{34}-t^1_4\delta^{34}_{13}+t^1_3\delta^{34}_{14})
\end{align*}}
\end{minipage}\\
\end{tabularx}
\end{sideways} 
\end{table}

%
%
%
%

\subsection{A non-parallelisable example}\label{nonexam}

The Kuranishi space of a nilmanifold which is  neither complex parallelisable nor carries an abelian complex structure can be much more complicated. We will illustrate this fact by describing a 2-step nilpotent Lie-algebra such that the Kuranishi space of an associated nilmanifold is singular but not cut out by quadrics, i.e., there are non-vanishing obstructions of higher order.

We use here an alternative way to describe a real Lie-algebra with complex structure: consider the complex vectorspace $V:=\langle X_1, \dots, X_7\rangle_\IC$. There is a natural real vectorspace $\lh\subset V\oplus \bar V$ invariant under complex conjugation such that $\lh_\IC=V\oplus \bar V$. This decomposition defines a complex structure $J$ on $\lh$ via  $\einsnull \lh:=V$.

Let $\omega^1, \dots, \omega^7$ be the basis of $V^*$ dual to the $X_i$. Then, by the formula for the differential \refb{differential}, a Lie bracket on $\lh$ is uniquely determinded by
\begin{gather*}
d\omega^1=d\omega^2=d\omega^3=d\omega^4=d\omega^5=0,\\
d\omega^6=\omega^1\wedge\omega^2,\\
d\omega^7=\omega^3\wedge\omega^4+ \bar \omega^1\wedge\omega^5,
\end{gather*}
and the complex conjugate equations.
For example, we have $[\bar X_5, X_1]=\bar X_7$.

Then 
\[d \einsnull {\lh^*}\subset \Lambda^{2,0}\oplus \Lambda^{1,1}\]
which means $d=\del+\delbar$ and the complex structure is integrable with respect to this Lie bracket. But since the image of $d$ is not contained in one of the components  $\Lambda^{1,1}$ and $\Lambda^{2,0}$ neither the complex structure is abelian nor is $(\lh, J)$ a complex Lie-algebra.

Our Lie-algebra with complex structure $(\lh,J)$ is defined over $\IQ$ and by the theorem of Mal'cev \cite{malcev51} there exists a lattice $\Gamma$ in the corresponding real simply connected nilpotent Lie-group $H$. We obtain a nilmanifold with left-invariant complex structure $(M,J)=(\Gamma\backslash H, J)$. 

Now let \[\mu:= \bar \omega^3\tensor X_1 +\bar \omega^4\tensor X_2.\]

Recall that for $X\in \einsnull\lh$ and $\bar Y \in \nulleins\lh$ we have $\delbar X (\bar Y)= \einsnull{ [\bar Y, X]}$
where $\einsnull x$ is the image of $x\in \lh_\IC$ under the projection to the $(1,0)$-part. In particular we see that $\delbar X_1=\delbar X_2=0$.
This implies $\delbar\mu=0$ and $\mu$ defines a class in $H^1((M,J), \Theta_{(M,J)})$.

Since every left-invariant function is constant and the contraction of a vector of type $(1,0)$ with a form of type $(0,2)$ is zero  the Schouten-bracket is given by
\[ [\bar\alpha\tensor X, \bar\beta\tensor Y]:=\bar\beta\wedge (i_Y\del\bar \alpha) \tensor X+ \bar\alpha \wedge (i_X\del\bar\beta)\tensor Y+\bar\alpha\wedge \bar \beta\tensor [X,Y].\]

We compute the first two steps of the iterative solution $\Phi$ with $\Phi_1=\mu$ of the Maurer-Cartan equation.

Since $\del\bar\omega^3=\del\bar\omega^4=0$ we get
\begin{align*}
[\mu,\mu]&= 2\bar\omega^3\wedge\bar\omega^4\tensor [X_1, X_2]\\
&=-\delbar(2\bar \omega^7\tensor X_6).
\end{align*}
We see that the obstruction in degree 2 vanishes.

Following the recursion \refb{Phi} we set $\Phi_2=2\bar\omega^7\tensor X_6$ and hence
 \begin{align*}
 [\Phi_1, \Phi_2]&= [\bar \omega^3\tensor X_1, 2\bar\omega^7\tensor X_6] +[\bar \omega^4\tensor X_2, 2\bar\omega^7\tensor X_6]\\
&= 2\bar\omega^3\wedge(i_{X_1}\del\bar\omega^7)\tensor X_6 +2\bar\omega^4\wedge(i_{X_2} \omega^1\wedge\bar\omega^5)\tensor X_6\\
&= 2\bar\omega^3\wedge\bar\omega^5\tensor X_6.
\end{align*}

It is immediate from the equations that this $2$-form with values in the tangent bundle is not exact and hence there is a non-vanishing obstruction in degree three.

%
%
%
%

\end{document}